\documentclass[runningheads]{llncs}

	\usepackage{latexsym}
	\usepackage{mathtools}
	\usepackage{amsmath}
	\usepackage{amssymb}
	\usepackage{stmaryrd}
	\usepackage{enumerate}
	\usepackage{color}
	\usepackage{graphicx}
	\usepackage{algorithm}
	\usepackage{url}
	\usepackage[colorlinks=true]{hyperref}
	\usepackage[numbers,sort&compress]{natbib}
	\usepackage{thmtools}
	\usepackage{tikz}

	\newcommand{\by}{\times}
	\newcommand{\intset}[1]{\left\llbracket #1 \right\rrbracket}

	\newcommand{\norm}[1]{\ensuremath{\left\lVert #1 \right\rVert}}
	\newcommand{\ip}[1]{\ensuremath{\left\langle #1 \right\rangle}}

	\newcommand{\R}{{\mathbb{R}}}

	\DeclareMathOperator{\Opt}{Opt}

	\DeclareMathOperator{\tr}{tr}
	\def\Se{{\mathbb{S}}}
	
	\DeclareMathOperator{\aff}{aff}

	\DeclareMathOperator{\conv}{conv}
	\DeclareMathOperator{\cone}{cone}
	\def\extr{{\mathop{\rm extr}}}

	\newcommand{\set}[1]{\left\{#1\right\}}

  \newcommand{\mb}{\mathbf}
  
  \newcommand{\mc}{\mathcal}
  \newcommand{\bb}{\mathbb}

\spnewtheorem{assumption}{Assumption}{\bf}{\it}
\spnewtheorem{observation}{Observation}{\bf}{\it}

\begin{document}

\title{On convex hulls of epigraphs of QCQPs
}
\titlerunning{Convex hulls of epigraphs of QCQPs
}

\author{
Alex L.\ Wang\inst{1}
\and
Fatma K{\i}l{\i}n\c{c}-Karzan\inst{1}
}
\authorrunning{A. Wang and F. K{\i}l{\i}n\c{c}-Karzan}

\institute{Carnegie Mellon University, Pittsburgh, PA, 15213, USA}

\maketitle

\begin{abstract}
Quadratically constrained quadratic programs (QCQPs) are a fundamental class of optimization problems well-known to be NP-hard in general. In this paper we study sufficient conditions for a \textit{convex hull result} that immediately implies that the standard semidefinite program (SDP) relaxation of a QCQP is tight.
We begin by outlining a general framework for proving such sufficient conditions. Then using this framework, we show that the convex hull result holds whenever the quadratic eigenvalue multiplicity, a parameter capturing the amount of symmetry present in a given problem, is large enough. 
Our results also imply new sufficient conditions for the tightness (as well as convex hull exactness) of a second order cone program relaxation of simultaneously diagonalizable QCQPs.

\end{abstract}

\section{Introduction}
In this paper we study \textit{quadratically constrained quadratic programs} (QCQPs) of the following form
\begin{align}
\label{eq:qcqp}
\Opt \coloneqq \inf_{x\in\R^N}\set{q_0(x) :\, \begin{array}
	{l}
q_i(x)\leq 0 ,\,\forall i\in\intset{m_I}\\
q_i(x) = 0 ,\,\forall i\in\intset{m_I+1,m_I+m_E}
\end{array}
},
\end{align}
where for every $i\in\intset{0,m_I+m_E}$, the function $q_i:\R^N\to\R$ is a (possibly nonconvex) quadratic function.
We will write $q_i(x) = x^\top A_i x + 2b_i^\top x + c_i$ where $A_i\in \bb S^N$, $b_i\in\R^N$, and $c_i\in\R$.
We will assume that the number of constraints $m\coloneqq m_I+m_E$ is at least $1$.

QCQPs arise naturally in many areas. A non-exhaustive list of applications contains facility location, production planning, pooling, max-cut, max-clique, and certain robust optimization problems (see \cite{phan1982quadratically,bao2011semidefinite,BenTal_ElGhaoui_Nemirovski_09} and references therein).

Although QCQPs are NP-hard to solve in general, they admit tractable convex relaxations. One natural relaxation is the standard (Shor) semidefinite program (SDP) relaxation \cite{shor1990dual}.
There is a vast literature on approximation guarantees associated with this relaxation \cite{BN2001,ye1999approximating,10.1007/978-3-0348-8362-7_15,nesterov1997quality}, however, less is known about its exactness.
Recently, a number of exciting results in phase retrieval~\cite{candes2015phase} and clustering~\cite{mixon2016clustering,abbe2015exact,rujeerapaiboon2019size} 
have shown that under various assumptions on the data,
the QCQP formulation of the corresponding problem has a tight SDP relaxation.
In contrast to these results, which address QCQPs arising from particular problems, \citet{burer2018exact}  very recently gave appealing deterministic sufficient conditions under which the standard SDP relaxation of \emph{general} QCQPs is tight. 
In our paper, 
we continue this vein of research for general QCQPs.   
More precisely, we will provide sufficient conditions under which the \emph{convex hull} of the epigraph of the QCQP is given by the
projection of the epigraph of its SDP relaxation. Note that such a result immediately implies that the optimal objective value of the QCQP is equal to the optimal objective value of its SDP relaxation.
We will refer to these two types of results as ``convex hull results'' and ``SDP tightness results.''
In this paper we will focus mainly on conditions that imply the convex hull result. See the full paper \cite{wang2019tightness} for additional new conditions which imply the SDP tightness result directly.

Convex hull results will necessarily require stronger assumptions than SDP tightness results, however they are also more broadly applicable because
they may be used
to derive strong convex relaxations for complex problems. 
In fact, the convexification of commonly occurring substructures has been critical in advancing the state-of-the-art computational approaches
for mixed integer linear programs and general nonlinear nonconvex programs~\cite{conforti2014integer,tawarmalani2002convexification}.
For computational purposes, conditions guaranteeing simple convex hull descriptions are particularly favorable.
As we will discuss later, a number of our sufficient conditions will guarantee that the desired convex hulls are given by a finite number of easily computable convex quadratic constraints in the original space of variables.

\subsubsection{Related work}

Convex hull results are well-known for simple QCQPs such as the Trust Region Subproblem (TRS) and the Generalized Trust Region Subproblem (GTRS).
Recall that the TRS is a QCQP with a single strictly convex inequality constraint and that the GTRS is a QCQP with a single (possibly nonconvex) inequality constraint.
A celebrated result due to \citet{FradkovYakubovich1979} implies that the SDP relaxation of the GTRS is tight.
More recently, \citet{Ho-NguyenKK2017} and \citet{wang2019generalized} showed that the (closed) convex hulls of the TRS and GTRS epigraphs are given exactly by the projection of the SDP epigraphs.
In both cases, the projections of the SDP epigraphs can also be described in the original space of variables with at most two convex quadratic inequalities.
As a result, the TRS and the GTRS can be solved without explicitly running costly SDP-based algorithms.

A different line of research has focused on providing explicit descriptions for the convex hull of the intersection of a single nonconvex quadratic region with convex sets such as convex quadratic regions, second-order cones~(SOCs), or polytopes, or with one other nonconvex quadratic region~\cite{BKK14,KKY15,YC15,yildiran2009convex,modaresi2017convex,santana2018convex}.
For example, the convex hull of the intersection of a two-term disjunction, which is a nonconvex quadratic constraint under mild assumptions, with the second-order cone (SOC) or its cross sections has received much attention in mixed integer programming; see~\cite{BKK14,KKY15,YC15} and references therein.
In contrast to these results, we will not limit the number of nonconvex quadratic constraints in our QCQPs. On the other hand, the nonconvex sets that we study in this paper will arise as epigraphs of QCQPs. In particular, the epigraph variable will play a special role in our analysis.
Therefore, we view our developments as complementary to these results.

The convex hull question has also received attention for certain strengthened relaxations of simple QCQPs \cite{SturmZhang2003,BurerAnstreicher2013,BurerYang2014,Burer2015}.
In this line of work, the standard SDP relaxation is strengthened by additional inequalities derived using the Reformulation-Linearization Technique (RLT). 
For example, \citet{SturmZhang2003} showed that the standard SDP relaxation strengthened with an additional SOC constraint derived from RLT gives the convex hull of the epigraph of the TRS with one additional linear inequality.
See \cite{Burer2015} for a survey of some results in this area.
In this paper, we restrict our attention to the standard SDP relaxation of QCQPs. Nevertheless, exactness conditions for strengthened SDP relaxations of QCQPs are clearly of great interest and are a direction for future research.

A number of SDP tightness results are known for variants of the TRS~\cite{JeyakumarLi2013,Ho-NguyenKK2017,YeZhang2003,BeckEldar2006}, for simultaneously diagonalizable QCQPs~\cite{Locatelli2016}, quadratic matrix programs~\cite{beck2007quadratic,beck2012new}, and random general QCQPs~\cite{burer2018exact}. See the full version of this paper for a more complete survey of the related SDP tightness results.

\subsubsection{Overview and outline of paper}

In contrast to the literature, which has mainly focused on simple QCQPs or QCQPs under certain structural assumptions, in this paper, we will consider general QCQPs and develop sufficient conditions for both the convex hull result and the SDP tightness result. 

We first introduce the epigraph of the QCQP by writing
\begin{align*}
\Opt &= \inf_{(x,t)\in\R^{N+1}} \set{2t :\, (x,t)\in\mc D},
\end{align*}
where $\mc D$ is the epigraph of the QCQP in \eqref{eq:qcqp}, i.e.,
\begin{align}
\label{eq:qcqp_epi}
\mc D \coloneqq \set{(x,t)\in\R^N\times \R :\, \begin{array}
	{l}
	q_0(x) \leq 2t\\
	q_i(x) \leq 0 ,\,\forall i\in\intset{m_I}\\
	q_i(x) = 0,\,\forall i\in\intset{m_I+1,m}
\end{array}}.
\end{align}

As $(x,t)\mapsto 2t$ is linear, we may replace the (potentially nonconvex) epigraph $\mc D$ with its convex hull $\conv(\mc D)$. Then,
\begin{align*}
\Opt &= \inf_{(x,t)\in\R^{N+1}} \set{2t :\, (x,t)\in\conv(\mc D)}.
\end{align*}

A summary of our contributions\footnote{
Due to space constraints, we omit full proofs, more detailed comparisons of our results with the literature, and our SDP tightness results in this extended abstract. The full version of this paper can be found at~\cite{wang2019tightness}.
}, along with an outline of the paper, is as follows. 
In Section~\ref{sec:framework}, we introduce and study the standard SDP relaxation of QCQPs \cite{shor1990dual} along with its optimal value $\Opt_\textup{SDP}$ and projected epigraph $\mc D_\textup{SDP}$. We set up a framework for deriving sufficient conditions for the ``convex hull result,'' $\conv(\mc D) = \mc D_\textup{SDP}$, and the ``SDP tightness result,'' $\Opt=\Opt_\textup{SDP}$. This framework is based on the Lagrangian function $(\gamma,x)\mapsto q_0(x) + \sum_{i=1}^m \gamma_i q_i(x)$ and the eigenvalue structure of a dual object $\Gamma\subseteq\R^m$.
This object $\Gamma$, which consists of the convex Lagrange multipliers, has been extensively studied in the literature (see \cite[Chapter 13.4]{wolkowicz2012handbook} and more recently \cite{sheriff2013convexity}).
In Section~\ref{sec:symmetries}, we define an integer parameter $k$, the quadratic eigenvalue multiplicity, that captures the amount of symmetry in a given QCQP. We then give examples where the quadratic eigenvalue multiplicity is large.
Specifically, vectorized reformulations of quadratic matrix programs \cite{beck2007quadratic} are such an example.
In Section~\ref{sec:conv_hull}, we use our framework to derive sufficient conditions for the convex hull result:~$\conv(\mc D) = \mc D_\textup{SDP}$. Theorem~\ref{thm:conv_hull_symmetries} states that if $\Gamma$ is polyhedral and $k$ is sufficiently large, then $\conv(\mc D) = \mc D_\textup{SDP}$. This theorem  actually follows as a consequence of Theorem~\ref{thm:conv_hull_main}, which replaces the assumption on the quadratic eigenvalue multiplicity with a weaker assumption regarding the dimension of zero eigenspaces related to the $A_i$ matrices.
Furthermore, our results in this section establish that if $\Gamma$ is polyhedral, then $\mc D_\textup{SDP}$ is SOC representable; see Remark~\ref{rem:soc_representability}. In particular, when the assumptions of Theorems~\ref{thm:conv_hull_main}~or~\ref{thm:conv_hull_symmetries} hold, we have that $\conv(\mc D) = \mc D_\textup{SDP}$ is SOC representable.
We provide several classes of problems that satisfy the assumptions of these theorems.
In particular, we recover a number of results regarding the TRS~\cite{Ho-NguyenKK2017}, the GTRS~\cite{wang2019generalized}, and the solvability of systems of quadratic equations \cite{barvinok1993feasibility}.

To the best of our knowledge, our results are the first to provide a unified explanation of many of the exactness guarantees in the literature.
Moreover, we provide significant generalizations of known results in a number of settings.

\subsubsection{Notation}
For nonnegative integers $m\leq n$ let $\intset{n}\coloneqq\{1,\ldots,n\}$ and $\intset{m,n}\coloneqq \set{m,m+1,\dots,n-1,n}$. 
Let $\mb{S}^{n-1}=\set{x\in\R^n:\, \norm{x}=1}$ denote the $n-1$ sphere.
Let $\Se^n$ denote the set of real symmetric $n\by n$ matrices.
For a positive integer $n$, let $I=I_n$ denote the $n\times n$ identity matrix. 
When the dimension is clear, we will simply write $I$.
Given two matrices $A$ and $B$, let $A \otimes B$ denote their Kronecker product.  
For a set $\mc D \subseteq \R^n$, let $\conv(\mc D)$, $\cone(\mc D)$, $\extr(\mc D)$,
$\dim(\mc D)$ and $\aff\dim(\mc D)$ denote the convex hull, conic hull, extreme points,
dimension, and affine dimension of $\mc D$, respectively.

\section{A general framework}
\label{sec:framework}

In this section, we introduce a general framework for analyzing the standard Shor SDP relaxation of QCQPs. We will examine how both the objective value and the feasible domain change when moving from a QCQP to its SDP relaxation.

We make an assumption that can be thought of as a primal feasibility and dual strict feasibility assumption. This assumption (or a slightly stronger version of it) is standard and is routinely made in the literature on QCQPs~\cite{BenTalTeboulle1996,YeZhang2003,beck2007quadratic}.
\begin{assumption}
\label{as:gamma_definite}
	Assume the feasible region of \eqref{eq:qcqp} is nonempty and there exists $\gamma^*\in\R^m$ such that $\gamma^*_i\geq 0$ for all $i\in\intset{m_I}$ and $A_0 + \sum_{i=1}^m \gamma^*_i A_i \succ 0$.
\qed\end{assumption}

The standard SDP relaxation to~\eqref{eq:qcqp} is
\begin{align}
\label{eq:shor_sdp}
\Opt_\textup{SDP} \coloneqq \inf_{x\in\R^N, X\in \Se^{N}} \set{\ip{Q_0, Y}:\, \begin{array}
	{l}
	Y\coloneqq \begin{pmatrix}
		1 & x^\top \\ x & X
	\end{pmatrix}\\
	\ip{Q_i, Y}\leq 0 ,\,\forall i\in\intset{m_I}\\
	\ip{Q_i, Y}= 0 ,\,\forall i\in\intset{m_I+1,m}\\
	Y\succeq 0
\end{array}},
\end{align}
where $Q_i\in\Se^{N+1}$ is the matrix $Q_i \coloneqq \left(\begin{smallmatrix}
	c_i & b_i^\top \\
	b_i & A_i
\end{smallmatrix}\right)$.
Let $\mc D_{\textup{SDP}}$ denote the epigraph of the relaxation \eqref{eq:shor_sdp} projected away from the $X$ variables, i.e., define
\begin{align}
\label{eq:sdp_epi}
\mc D_{\textup{SDP}} \coloneqq \set{(x,t) \in\R^{N+1}:\, \begin{array}
	{l}
	\exists X\in \Se^N:\\
	Y \coloneqq\begin{pmatrix}
		1 & x^\top\\ x & X
	\end{pmatrix}\\
	\ip{Q_0, Y} \leq 2t\\
	\ip{Q_i, Y}\leq 0 ,\,\forall i\in\intset{m_I}\\
	\ip{Q_i, Y}= 0 ,\,\forall i\in\intset{m_I+1,m}\\
	Y\succeq 0
\end{array}}.
\end{align}

By taking $X = xx^\top$ in both \eqref{eq:shor_sdp} and \eqref{eq:sdp_epi}, we see that $\mc D \subseteq \mc D_\textup{SDP}$ and $\Opt\geq \Opt_\textup{SDP}$. Noting that $\mc D_\textup{SDP}$ is convex (it is the projection of a convex set), we further have that $\conv(\mc D)\subseteq \mc D_{\textup{SDP}}$. The framework that we set up in the remainder of this section allows us to reason about when equality occurs in either relation.

\subsection{Rewriting the SDP in terms of a dual object}

For $\gamma\in\R^m$, define
\begin{align*}
&A(\gamma)\coloneqq A_0 + \sum_{i=1}^m \gamma_i A_i,
\quad
b(\gamma)\coloneqq b_0 + \sum_{i=1}^m \gamma_i b_i,
\quad
c(\gamma)\coloneqq c_0 + \sum_{i=1}^m \gamma_i c_i,
\\
&q(\gamma,x)\coloneqq q_0(x) + \sum_{i=1}^m \gamma_i q_i(x).
\end{align*}

Our framework for analyzing \eqref{eq:shor_sdp} is based on the dual object
\begin{align*}
\Gamma\coloneqq \set{\gamma\in\R^m:\, \begin{array}
	{l}
	A(\gamma)\succeq 0\\
	\gamma_i\geq 0 ,\,\forall i\in\intset{m_I}
\end{array}}.
\end{align*}
This object will play a key role our analysis for the following fundamental reason.

\begin{lemma}
\label{lemma:sdp_in_terms_of_Gamma}
Suppose Assumption~\ref{as:gamma_definite} holds. Then
\begin{align*}
\mc D_\textup{SDP}  = \set{(x,t) :~ \sup_{\gamma\in\Gamma} q(\gamma,x)\leq 2t}
\quad\text{and}\quad
\Opt_\textup{SDP} &= \min_{x\in\R^N} \sup_{\gamma\in\Gamma} q(\gamma, x).
\end{align*}
\end{lemma}
The second identity is well-known; see e.g., \citet{Fujie1997}.

\subsection{The eigenvalue structure of $\Gamma$}
We now define a number of objects related to $\Gamma$. Noting that $\gamma\mapsto q(\gamma,\hat x)$ is linear and that $\Gamma$ is closed leads to the following observation.

\begin{observation}\label{obs:GammaFace}
Let $\hat x\in\R^N$. If $\sup_{\gamma\in\Gamma} q(\gamma,\hat x)$ is finite, then $q(\gamma,\hat x)$ achieves its maximum value in $\Gamma$ on some face $\mc F$ of $\Gamma$.
\end{observation}
In particular, the following definition is well-defined.

\begin{definition}
For any $\hat x\in\R^N$ such that $\sup_{\gamma\in\Gamma} q(\gamma,\hat x)$ is finite, define $\mc F(\hat x)$ to be the face of $\Gamma$ maximizing $q(\gamma,\hat x)$.
\end{definition}

\begin{definition}\label{def:definiteFace}
Let $\mc F$ be a face of $\Gamma$. We say that $\mc F$ is a \textit{definite face} if there exists $\gamma\in \mc F$ such that $A(\gamma)\succ 0$. Otherwise, we say that $\mc F$ is a \textit{semidefinite face} and let $\mc V(\mc F)$ denote the shared zero eigenspace of $\mc F$, i.e.,
\begin{align*}
\mc V(\mc F) &\coloneqq \set{v\in\R^N:\, A(\gamma)v = 0,\,\forall \gamma\in\mc F}.
\end{align*}
\end{definition}

It is possible to show that for $\mc F$ semidefinite, the set $\mc V(\mc F)$ is nontrivial. As a sketch, suppose otherwise, then for every $v$ on the unit sphere, we can associate a $\gamma_v\in \mc F$ such that $v^\top A(\gamma_v)v >0$. Then we can produce a positive definite matrix $A(\bar\gamma)$ where $\bar\gamma$ is an ``average'' over the $\gamma_v$, a contradiction. See Lemma 2 in the full version of this paper for a formal proof.

\subsection{The framework}
Our framework
consists of two parts: an ``easy part'' that only requires Assumption~\ref{as:gamma_definite} to hold and a ``hard part'' that may require much stronger assumptions. The ``easy part'' consists of the following lemma and observation.

\begin{lemma}
\label{lem:F_contains_definite}
Suppose Assumption~\ref{as:gamma_definite} holds and let $(\hat x, \hat t) \in\mc D_\textup{SDP}$. If $\mc F(\hat x)$ is a definite face of $\Gamma$, then $(\hat x, \hat t)\in\mc D$.
\end{lemma}

\begin{observation}
\label{obs:aff_dim_m_is_definite}
Suppose Assumption~\ref{as:gamma_definite} holds and let $\mc F$ be a face of $\Gamma$. If $\aff\dim(\mc F) = m$, then $\mc F$ is definite.
\end{observation}

The ``hard part'' of the framework works as follows:
In order to show the convex hull result $\mc D_{\textup{SDP}}= \conv(\mc D)$, it suffices to guarantee that every $(\hat x,\hat t)\in\mc D_\textup{SDP}$ can be decomposed as a convex combination of pairs $(x_\alpha,t_\alpha)$ for which $\mc F(x_\alpha)$ is definite. Then, by Lemma~\ref{lem:F_contains_definite}, we will have that $(x_\alpha,t_\alpha)\in\mc D$.
We give examples of such sufficient conditions in Section~\ref{sec:conv_hull}. Our decomposition procedures will be recursive and we will use Observation~\ref{obs:aff_dim_m_is_definite} to show that they terminate.

\begin{remark}
\label{rem:affine_transformation}
Consider performing an invertible affine transformation on the space $\R^N$, i.e. let $y = U(x + z)$ where $U\in\R^{N\by N}$ is an invertible linear transformation and $z\in\R^N$. Define the quadratic functions $q'_0,\dots,q'_m: \R^N\to\R$ such that $q'_i(y) = q'_i(U(x+z))=q_i(x)$ for all $x\in\R^N$. 
We will use an apostrophe to denote all the quantities corresponding to the QCQP in the variable $y$.

Define the map $\ell:\R^{N+1}\to\R^{N+1}$ by $(x,t)\mapsto (U(x+z), t)$. Note that $\Opt' = \Opt$ and $\conv(\mc D') = \ell(\conv(\mc D))$.
Furthermore a straightforward application of Lemma~\ref{lemma:sdp_in_terms_of_Gamma} gives
$\Opt_\textup{SDP}' = \Opt_\textup{SDP}$ and $\mc D_\textup{SDP}' = \ell(\mc D_\textup{SDP})$. We deduce that the questions $\conv(\mc D)\stackrel{?}{=} \mc D_\textup{SDP}$ and $\Opt\stackrel{?}{=} \Opt_\textup{SDP}$ are invariant under invertible affine transformation of the $x$-space.
In particular, the sufficient conditions that we will present in Theorems~\ref{thm:conv_hull_main}~and~\ref{thm:conv_hull_symmetries} only need to hold after some invertible affine transformation. In this sense, the SDP relaxation will ``find'' structure in a given QCQP even if it is ``hidden'' by an affine transformation.
\qed\end{remark}

\section{Symmetries in QCQPs}
\label{sec:symmetries}

In this section, we examine a parameter $k$ that captures the amount of symmetry present in a QCQP of the form \eqref{eq:qcqp}.

\begin{definition}
The \textit{quadratic eigenvalue multiplicity} of a QCQP of the form \eqref{eq:qcqp} is the largest integer $k$ such that for every $i\in\intset{0,m}$ there exists $\mc A_i \in \bb S^n$ for which $A_i = I_k\otimes \mc A_i$.
Let $\mc A(\gamma) \coloneqq \mc A_0 + \sum_{i=1}^m \gamma_i \mc A_i$.
\end{definition}
This value is well-defined: $k$ is always at least $1$ as we can write $A_i = I_1 \otimes \mc A_i$. On the other hand, $k$ must also be a divisor of $N$.

The next lemma states the crucial structure inherent in QCQPs with large quadratic eigenvalue multiplicities.
\begin{lemma}
\label{lem:V_F_large}
If $\mc F$ is a semidefinite face of $\Gamma$, then $\dim(\mc V(\mc F)) \geq k$.
\end{lemma}

\begin{remark}
In quadratic matrix programming \cite{beck2007quadratic,beck2012new}, we are asked to optimize
\begin{align}
\label{eq:quadratic_matrix_programming}
\inf_{X\in\R^{n\by k}}\set{
\begin{array}
	{l}
	\tr(X^\top \mc A_0 X)+2\tr(B_0^\top X) + c_0:\,\\
	\qquad \tr(X^\top \mc A_i X)+2\tr(B_i^\top X) + c_i \leq 0,\,\forall i\in\intset{m_I}\\
	\qquad \tr(X^\top \mc A_i X)+2\tr(B_i^\top X) + c_i = 0,\,\forall i\in\intset{m_I+1,m}
\end{array}},
\end{align}
where $\mc A_i\in\bb S^n$, $B_i\in\R^{n\by k}$ and $c_i\in\R$ for all $i\in\intset{0,m}$. We can transform this program to an equivalent QCQP in the vector variable $x\in\R^{nk}$.
Then
$\tr(X^\top \mc A_i X)+2\tr(B_i^\top X) + c_i = x^\top \left(I_k\otimes \mc A_i\right)x+ 2 b_i^\top x + c_i$,
where $b_i\in\R^{nk}$ has entries $(b_i)_{(t-1)n + s} = (B_i)_{s,t}$. In particular, the vectorized reformulation of \eqref{eq:quadratic_matrix_programming} has quadratic eigenvalue value multiplicity $k$.
\qed\end{remark}

\section{Convex hull results}
\label{sec:conv_hull}
We now present new sufficient conditions for the convex hull result $\mc D_\textup{SDP} = \conv(\mc D)$. We analyze the case where the geometry of $\Gamma$ is particularly nice.

\begin{assumption}
\label{as:gamma_polyhedral}
	Assume that $\Gamma$ is polyhedral.
\qed\end{assumption}
We remark that although Assumption~\ref{as:gamma_polyhedral} is rather restrictive, it is general enough to cover the case where the set of quadratic forms $\set{A_i}_{i\in\intset{0,m}}$ is diagonal or simultaneously diagonalizable --- a class of QCQPs which have been studied extensively in the literature~\cite{BenTalDenHertog2014,locatelli2015some,Locatelli2016}.
See the full version of this paper for convex hull and SDP tightness results without Assumption~\ref{as:gamma_polyhedral} as well as a discussion on the difficulties in removing it.

Our main result in this paper is the following theorem.
\begin{theorem}
\label{thm:conv_hull_main}
Suppose Assumptions~\ref{as:gamma_definite}~and~\ref{as:gamma_polyhedral} hold. If for every semidefinite face $\mc F$ of $\Gamma$ we have
$\dim(\mc V(\mc F)) \geq \aff\dim(\set{b(\gamma):\,\gamma\in\mc F}) + 1$,
then
$\conv(\mc D) = \mc D_\textup{SDP}$.
\end{theorem}

Assumption~\ref{as:gamma_definite} allows us to apply Lemma~\ref{lem:F_contains_definite} to handle any $(\hat x,\hat t)\in\mc D_\textup{SDP}$ for which $\mc F(\hat x)$ is definite. Therefore, in order to prove Theorem~\ref{thm:conv_hull_main}, it suffices to prove the following lemma.
\begin{lemma}
\label{lem:pivoting_F}
Suppose Assumptions~\ref{as:gamma_definite}~and~\ref{as:gamma_polyhedral} hold. Let $(\hat x,\hat t)\in\mc D_\textup{SDP}$ and let $\mc F = \mc F(\hat x)$. If $\mc F$ is a semidefinite face of $\Gamma$ and
$\dim(\mc V(\mc F)) \geq \aff \dim(\set{b(\gamma):\, \gamma\in\mc F})+1$,
then $(\hat x,\hat t)$ can be written as a convex combination of points $(x_\alpha,t_\alpha)$ satisfying the following properties:
\begin{enumerate}
	\item $(x_\alpha,t_\alpha)\in\mc D_\textup{SDP}$, and
	\item $\aff\dim(\mc F(x_\alpha))>\aff\dim(\mc F(\hat x ))$.
\end{enumerate}
\end{lemma}
We give a proof sketch of Lemma~\ref{lem:pivoting_F} in Appendix~\ref{app:proof_sketch_lem_pivoting}.

The proof of Theorem~\ref{thm:conv_hull_main} follows at once from Lemmas~\ref{lem:F_contains_definite}~and~\ref{lem:pivoting_F} and Observation~\ref{obs:aff_dim_m_is_definite}. Indeed, Lemma~\ref{lem:pivoting_F} guarantees that $\aff\dim(\mc F(x_\alpha))>\aff\dim(\mc F(\hat x))$. Thus, by Observation~\ref{obs:aff_dim_m_is_definite}, we will have successfully decomposed $(\hat x,\hat t)$ as a convex combination of $(x_\alpha,t_\alpha)$, where $(x_\alpha,t_\alpha)\in\mc D_\textup{SDP}$ and $\mc F(x_\alpha)$ is definite, after at most $m-1$ rounds of applying Lemma~\ref{lem:pivoting_F}. Finally, Lemma~\ref{lem:F_contains_definite} guarantees that each pair $(x_\alpha,t_\alpha)$ is an element of $\mc D$, the epigraph of the QCQP.

The next theorem follows as a corollary to Theorem~\ref{thm:conv_hull_main}.

\begin{theorem}
\label{thm:conv_hull_symmetries}
Suppose Assumptions~\ref{as:gamma_definite}~and~\ref{as:gamma_polyhedral} hold. If for every semidefinite face $\mc F$ of $\Gamma$ we have
$k\geq \aff\dim(\set{b(\gamma):\, \gamma\in\mc F}) + 1$,
then
$\conv(\mc D) = \mc D_\textup{SDP}$.
\end{theorem}

\begin{remark}\label{rem:soc_representability}
We remark that when $\Gamma$ is polyhedral (Assumption~\ref{as:gamma_polyhedral}), the set $\mc D_\textup{SDP}$ is actually SOC representable:
By the Minkowski-Weyl Theorem, we can decompose $\Gamma = \Gamma_e + \cone(\Gamma_r)$ where both $\Gamma_e$ and $\Gamma_r$ are polytopes. Let $\breve q(\gamma,x) = \sum_{i=1}^m \gamma_i q_i(x)$. Then, by Lemma~\ref{lemma:sdp_in_terms_of_Gamma} we can write
\begin{align*}
\mc D_\textup{SDP} = \set{(x,t):\, \sup_{\gamma\in\Gamma} q(\gamma,x)\leq 2t}
=  \set{(x,t):\,\begin{array}
	{l}
	q(\gamma_e,x) \leq 2t ,\,\forall \gamma_e\in\extr(\Gamma_e)\\
	\breve q(\gamma_f,x)\leq 0 ,\,\forall \gamma_f\in\extr(\Gamma_r)
\end{array}}.
\end{align*}
That is, $\mc D_\textup{SDP}$ is defined by finitely many convex quadratic inequalities. Thus the assumptions of Theorems~\ref{thm:conv_hull_main}~and~\ref{thm:conv_hull_symmetries} imply that $\conv(\mc D)$ is SOC representable.
\qed\end{remark}

We now give examples of problems where our assumptions hold.

\begin{corollary}\label{cor:m=1}
Suppose $m = 1$ and Assumption~\ref{as:gamma_definite} holds. Then,
$\conv(\mc D) = \mc D_\textup{SDP}$.
\end{corollary}
Corollary~\ref{cor:m=1} recovers results associated with the epigraph of the TRS\footnote{
Corollary~\ref{cor:m=1} fails to recover the full extent of \cite[Theorem 13]{Ho-NguyenKK2017}. Indeed, \cite[Theorem 13]{Ho-NguyenKK2017} also gives a description of the convex hull of the epigraph of the TRS with an additional conic constraint under some assumptions.
}
and the GTRS (see \cite[Theorem 13]{Ho-NguyenKK2017} and \cite[Theorems 1 and 2]{wang2019generalized}).

\begin{corollary}\label{cor:b_i=0}
Suppose Assumptions~\ref{as:gamma_definite}~and~\ref{as:gamma_polyhedral} hold. If $b_i= 0$ for all $i\in\intset{m}$, then
$\conv(\mc D) = \mc D_\textup{SDP}$.
\end{corollary}

\begin{example}
\label{ex:b_i=0}
Consider the following optimization problem.
\begin{align*}
\inf_{x\in\R^2} \set{x_1^2 +x_2^2 +10x_1 :\, \begin{array}
	{l}
	x_1^2 - x_2^2 -5 \leq 0\\
	-x_1^2 + x_2^2 -50 \leq 0
\end{array}}
\end{align*}
We check that the conditions of Corollary~\ref{cor:b_i=0} hold. Assumption~\ref{as:gamma_definite} holds as $A(0) = A_0 = I\succ 0$ and $x= 0$ is feasible. Next, Assumption~\ref{as:gamma_polyhedral} holds as
\begin{align*}
\Gamma &= \set{\gamma\in\R^2:\, \begin{array}
	{l}
	1+\gamma_1-\gamma_2\geq 0\\
	1-\gamma_1+\gamma_2\geq 0\\
	\gamma\geq 0
\end{array}}.
\end{align*}
One can verify that 
$\Gamma = \conv\left(\set{(0,0), (1,0), (0,1)}\right) + \cone(\set{(1,1)})$. 
Finally, we note that $b_1 =b_2=0$. Hence, Corollary~\ref{cor:b_i=0} and Remark~\ref{rem:soc_representability} imply that
\begin{align*}
\conv(\mc D) = \mc D_\textup{SDP} = \set{(x,t) :\, \begin{array}
	{l}
	x_1^2+x_2^2 +10 x_1\leq 2t\\
	2x_1^2 +10x_1 -5 \leq 2t\\
	2x_2^2 +10x_1 - 50 \leq 2t
\end{array}}.
\end{align*}
We plot $\mc D$ and $\conv(\mc D)=\mc D_\textup{SDP}$ in Figure~\ref{fig:example_convex_hull}.
\begin{figure}
  \centering
    \includegraphics[width=0.4\textwidth]{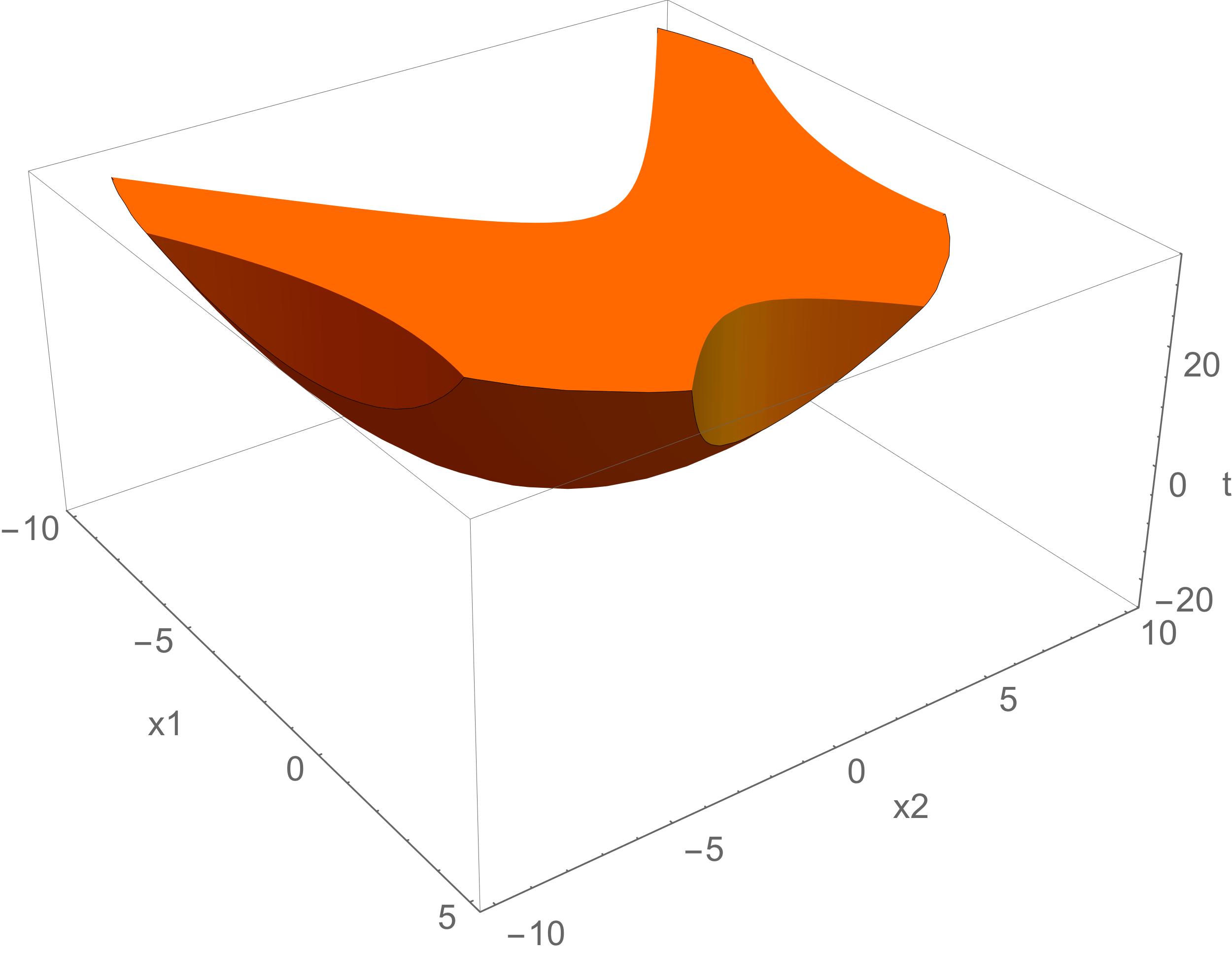}\qquad
    \includegraphics[width=0.4\textwidth]{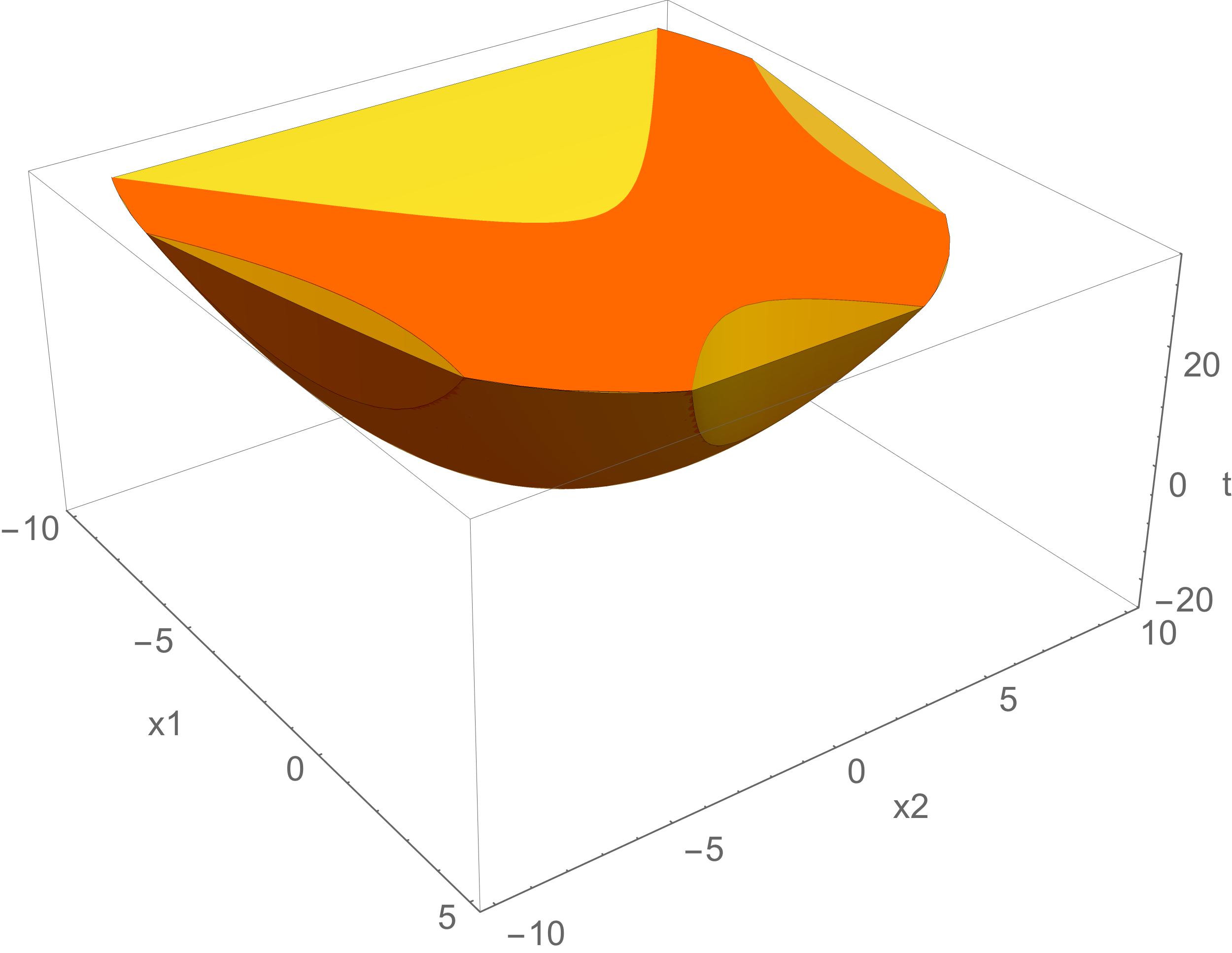}
	\caption{The sets $\mc D$ (in orange) and $\conv(\mc D)$ (in yellow) from Example~\ref{ex:b_i=0}}
	\label{fig:example_convex_hull}
\end{figure}
\qed\end{example}

\begin{remark}\label{ex:Barvinok}
\citet{barvinok1993feasibility} shows that one can decide in polynomial time (in $N$) whether a constant number, $m_E$, of quadratic forms $\set{A_i}_{i\in\intset{m_E}}$ has a joint nontrivial zero. That is, whether the system $x^\top A_i x=0$ for $i\in\intset{m_E}$ and $x^\top x = 1$ is feasible. We can recast this as asking whether the following optimization problem
\begin{align*}
\min_{x\in\R^N}\set{-x^\top x:~  \begin{array}
	{l}
	x^\top x\leq 1\\
	x^\top A_i x=0,\forall i\in\intset{m_E}
\end{array}}
\end{align*}
has objective value $-1$ or $0$.

Thus, the feasibility problem studied in \cite{barvinok1993feasibility} reduces to a QCQP of the form we study in this paper.
It is easy to verify that Assumption~\ref{as:gamma_definite} holds.
Then when $\Gamma$ is polyhedral (Assumption~\ref{as:gamma_polyhedral}), Corollary~\ref{cor:b_i=0} implies that the feasibility problem (even in a variable number of quadratic forms) can be decided using a semidefinite programming approach.
Nevertheless, Assumption~\ref{as:gamma_polyhedral} may not necessarily hold in general and so Corollary~\ref{cor:b_i=0}  does not recover the full result of \cite{barvinok1993feasibility}.
\qed\end{remark}

\begin{corollary}
\label{cor:identity_matrices}
Suppose Assumption~\ref{as:gamma_definite} holds and for every $i\in\intset{0,m}$, there exists $\alpha_i$ such that $A_i = \alpha_i I_N$. If $m\leq N$, then
$\conv(\mc D) = \mc D_\textup{SDP}$.
\end{corollary}

\begin{remark}
Consider the problem of finding the distance between the origin $0\in\R^N$ and a piece of Swiss cheese $C\subseteq \R^N$. We will assume that $C$ is nonempty and defined as
\begin{align*}
C = \set{x\in\R^N :\, \begin{array}
	{l}
	\norm{x - y_i} \leq s_i ,\,\forall i\in\intset{m_1}\\
	\norm{x - z_i} \geq t_i ,\,\forall i\in\intset{m_2}\\
	\ip{x,b_i} \geq c_i ,\,\forall i\in\intset{m_3}
\end{array}},
\end{align*}
where $y_i, z_i, b_i\in\R^N$ and $s_i,t_i,c_i\in\R$ are arbitrary.
In other words, $C$ is defined by $m_1$-many ``inside-ball'' constraints, $m_2$-many ``outside-ball'' constraints, and $m_3$-many linear inequalities. Note that each of these constraints may be written as a quadratic inequality constraint with quadratic form $I$, $-I$, or $0$. In particular, Corollary~\ref{cor:identity_matrices} implies that if $m_1+m_2+m_3\leq N$, then the value
\begin{align*}
\inf_{x\in\R^N}\set{\norm{x}^2 :\, x\in C}
\end{align*}
may be computed using the standard SDP relaxation of the problem.

\citet{BienstockMichalka2014} give sufficient conditions under which a related problem
\begin{align*}
\inf_{x\in\R^N}\set{q_0(x) :\, x\in C},
\end{align*}
is polynomial-time solvable. Here, $q_0:\R^N\to\R$ may be an arbitrary quadratic function however $m_1$ and $m_2$ must be constant.
They devise an enumerative algorithm for this problem and prove its correctness under different assumptions.
In contrast, our work deals only with the standard SDP relaxation and does not assume that the number of quadratic forms is constant.
\qed\end{remark}

\section*{Acknowledgments}
This research is supported in part by NSF grant CMMI 1454548.

\newpage
\appendix
\section{Proof sketch of Lemma~\ref{lem:pivoting_F}}
\label{app:proof_sketch_lem_pivoting}

For simplicity, we will assume that $\Gamma$ is a polytope in this proof sketch.
Let $(\hat x,\hat t)$ satisfy the assumptions of Lemma~\ref{lem:pivoting_F}. Without loss of generality, we may assume that $\sup_{\gamma\in\Gamma} q(\gamma,\hat x) = 2\hat t$.

We claim that the following system in variables $v$ and $s$
\begin{align*}
	\begin{cases}
		\ip{b(\gamma), v} = s,\,\forall \gamma\in \mc F\\
		v\in\mc V(\mc F),\,s\in\R
	\end{cases}
\end{align*}
has a nonzero solution. Indeed, we may replace the first constraint with at most
\begin{align*}
\aff\dim(\set{b(\gamma):\, \gamma\in\mc F}) + 1 \leq \dim(\mc V(\mc F))
\end{align*}
homogeneous linear equalities in the variables $v$ and $s$.
The claim then follows by noting that the equivalent system is an under-constrained homogeneous system of linear equalities and thus has a nonzero solution $(v,s)$. It is easy to verify that $v\neq 0$ and hence, by scaling, we may take $v\in\mb S^{N-1}$.

We will modify $(\hat x,\hat t)$ in the
$(v, s)$
direction. For $\alpha\in\R$, define
\begin{align*}
(x_\alpha,t_\alpha)\coloneqq \left(\hat x + \alpha v,\, \hat t + \alpha s\right).
\end{align*}
We will sketch the existence of an $\alpha>0$ such that $(x_\alpha,t_\alpha)$ satisfies the conclusions of Lemma~\ref{lem:pivoting_F}. A similar line of reasoning will produce an analogous $\alpha<0$. This will complete the proof sketch.

Suppose $\gamma\in\mc F$. Then, by our choice of $v$ and $s$, the function $\alpha\mapsto q(\gamma,x_\alpha)-2t_\alpha = q(\gamma,\hat x) -2t = 0$ is identically zero.
Now suppose $\gamma\in\Gamma\setminus\mc F$. Then, the function $\alpha\mapsto q(\gamma,x_\alpha)-2t_\alpha$ is a convex quadratic function which is negative at $\alpha = 0$.

We conclude that the following set
\begin{align*}
\mc Q \coloneqq \set{\alpha\mapsto q(\gamma,x_\alpha) - 2 t_\alpha :\, \gamma\in\extr(\Gamma)} \setminus \set{0},
\end{align*}
consists of convex quadratic functions which are negative at $\alpha = 0$. The finiteness of this set follows from the assumption that $\Gamma$ is polyhedral.

Assumption~\ref{as:gamma_definite} implies that at least one of the functions in $\mc Q$ is strictly convex. Then  as $\mc Q$ is a finite set,
there exists an $\alpha_+ >0$ such that $q(\alpha_+)\leq 0$ for all $q\in\mc Q$ with at least one equality. We emphasize that this is the step where Assumption~\ref{as:gamma_polyhedral} cannot be dropped.

Finally, it is easy to check that $(x_{\alpha_+}, t_{\alpha_+})$ satisfies the conclusions of Lemma~\ref{lem:pivoting_F}.

\bibliography{repeated_eig_bib.bib}

\begin{thebibliography}{43}
\providecommand{\natexlab}[1]{#1}
\providecommand{\url}[1]{\texttt{#1}}
\providecommand{\urlprefix}{URL }
\expandafter\ifx\csname urlstyle\endcsname\relax
  \providecommand{\doi}[1]{doi:\discretionary{}{}{}#1}\else
  \providecommand{\doi}{doi:\discretionary{}{}{}\begingroup
  \urlstyle{rm}\Url}\fi

\bibitem[{Abbe et~al.(2015)Abbe, Bandeira, and Hall}]{abbe2015exact}
Abbe, E., Bandeira, A.S., Hall, G.: Exact recovery in the stochastic block
  model. IEEE Transactions on Information Theory \textbf{62}(1), 471--487
  (2015)

\bibitem[{Bao et~al.(2011)Bao, Sahinidis, and
  Tawarmalani}]{bao2011semidefinite}
Bao, X., Sahinidis, N.V., Tawarmalani, M.: Semidefinite relaxations for
  quadratically constrained quadratic programming: A review and comparisons.
  Mathematical programming \textbf{129}(1), 129 (2011)

\bibitem[{Barvinok(1993)}]{barvinok1993feasibility}
Barvinok, A.I.: Feasibility testing for systems of real quadratic equations.
  Discrete \& Computational Geometry \textbf{10}(1), 1--13 (1993)

\bibitem[{Beck(2007)}]{beck2007quadratic}
Beck, A.: Quadratic matrix programming. SIAM Journal on Optimization
  \textbf{17}(4), 1224--1238 (2007)

\bibitem[{Beck et~al.(2012)Beck, Drori, and Teboulle}]{beck2012new}
Beck, A., Drori, Y., Teboulle, M.: A new semidefinite programming relaxation
  scheme for a class of quadratic matrix problems. Operations Research Letters
  \textbf{40}(4), 298--302 (2012)

\bibitem[{Beck and Eldar(2006)}]{BeckEldar2006}
Beck, A., Eldar, Y.C.: Strong duality in nonconvex quadratic optimization with
  two quadratic constraints. SIAM Journal on Optimization \textbf{17}(3),
  844--860 (2006)

\bibitem[{Ben-Tal et~al.(2009)Ben-Tal, {El Ghaoui}, and
  Nemirovski}]{BenTal_ElGhaoui_Nemirovski_09}
Ben-Tal, A., {El Ghaoui}, L., Nemirovski, A.: Robust Optimization. Princeton
  University Press, Princeton Series in Applied Mathematics, Philadelphia, PA,
  USA (2009)

\bibitem[{Ben-Tal and den Hertog(2014)}]{BenTalDenHertog2014}
Ben-Tal, A., den Hertog, D.: Hidden conic quadratic representation of some
  nonconvex quadratic optimization problems. Mathematical Programming
  \textbf{143}(1), 1--29 (2014)

\bibitem[{Ben-Tal and Nemirovski(2001)}]{BN2001}
Ben-Tal, A., Nemirovski, A.: Lectures on Modern Convex Optimization. MPS-SIAM
  Series on Optimization, SIAM, Philadehia, PA, USA (2001)

\bibitem[{Ben-Tal and Teboulle(1996)}]{BenTalTeboulle1996}
Ben-Tal, A., Teboulle, M.: Hidden convexity in some nonconvex quadratically
  constrained quadratic programming. Mathematical Programming \textbf{72}(1),
  51–63 (1996)

\bibitem[{Bienstock and Michalka(2014)}]{BienstockMichalka2014}
Bienstock, D., Michalka, A.: Polynomial solvability of variants of the
  trust-region subproblem. In: Proceedings of the Twenty-Fifth Annual ACM-SIAM
  Symposium on Discrete Algorithms, pp. 380--390 (2014)

\bibitem[{Burer(2015)}]{Burer2015}
Burer, S.: A gentle, geometric introduction to copositive optimization.
  Mathematical Programming \textbf{151}(1), 89--116 (2015)

\bibitem[{Burer and Anstreicher(2013)}]{BurerAnstreicher2013}
Burer, S., Anstreicher, K.M.: Second-order-cone constraints for extended
  trust-region subproblems. SIAM Journal on Optimization \textbf{23}(1),
  432--451 (2013)

\bibitem[{Burer and K{\i}l{\i}n\c{c}-Karzan(2017)}]{BKK14}
Burer, S., K{\i}l{\i}n\c{c}-Karzan, F.: How to convexify the intersection of a
  second order cone and a nonconvex quadratic. Mathematical Programming
  \textbf{162}(1), 393--429 (2017)

\bibitem[{Burer and Yang(2015)}]{BurerYang2014}
Burer, S., Yang, B.: The {Trust Region Subproblem} with non-intersecting linear
  constraints. Mathematical Programming \textbf{149}(1), 253--264 (2015)

\bibitem[{Burer and Ye(2018)}]{burer2018exact}
Burer, S., Ye, Y.: Exact semidefinite formulations for a class of (random and
  non-random) nonconvex quadratic programs. Mathematical Programming pp. 1--17
  (2018)

\bibitem[{Candes et~al.(2015)Candes, Eldar, Strohmer, and
  Voroninski}]{candes2015phase}
Candes, E.J., Eldar, Y.C., Strohmer, T., Voroninski, V.: Phase retrieval via
  matrix completion. SIAM review \textbf{57}(2), 225--251 (2015)

\bibitem[{Conforti et~al.(2014)Conforti, Cornu{\'e}jols, and
  Zambelli}]{conforti2014integer}
Conforti, M., Cornu{\'e}jols, G., Zambelli, G.: Integer programming, vol. 271.
  Springer (2014)

\bibitem[{Fradkov and Yakubovich(1979)}]{FradkovYakubovich1979}
Fradkov, A.L., Yakubovich, V.A.: The {S}-procedure and duality relations in
  nonconvex problems of quadratic programming. Vestn. LGU, Ser. Mat., Mekh.,
  Astron \textbf{6}(1), 101--109 (1979)

\bibitem[{Fujie and Kojima(1997)}]{Fujie1997}
Fujie, T., Kojima, M.: Semidefinite programming relaxation for nonconvex
  quadratic programs. Journal of Global Optimization \textbf{10}(4), 367--380
  (Jun 1997), ISSN 1573-2916, \doi{10.1023/A:1008282830093},
  \urlprefix\url{https://doi.org/10.1023/A:1008282830093}

\bibitem[{Phan-huy Hao(1982)}]{phan1982quadratically}
Phan-huy Hao, E.: Quadratically constrained quadratic programming: Some
  applications and a method for solution. Zeitschrift f{\"u}r Operations
  Research \textbf{26}(1), 105--119 (1982)

\bibitem[{Ho-Nguyen and K{\i}l{\i}n\c{c}-Karzan(2017)}]{Ho-NguyenKK2017}
Ho-Nguyen, N., K{\i}l{\i}n\c{c}-Karzan, F.: A second-order cone based approach
  for solving the {Trust Region Subproblem} and its variants. SIAM Journal on
  Optimization \textbf{27}(3), 1485--1512 (2017)

\bibitem[{Jeyakumar and Li(2014)}]{JeyakumarLi2013}
Jeyakumar, V., Li, G.Y.: Trust-region problems with linear inequality
  constraints: {E}xact {SDP} relaxation, global optimality and robust
  optimization. Mathematical Programming \textbf{147}(1), 171--206 (2014)

\bibitem[{K{\i}l{\i}n\c{c}-Karzan and Y{\i}ld{\i}z(2015)}]{KKY15}
K{\i}l{\i}n\c{c}-Karzan, F., Y{\i}ld{\i}z, S.: Two-term disjunctions on the
  second-order cone. Mathematical Programming \textbf{154}(1), 463--491 (2015)

\bibitem[{Locatelli(2015)}]{locatelli2015some}
Locatelli, M.: Some results for quadratic problems with one or two quadratic
  constraints. Operations Research Letters \textbf{43}(2), 126--131 (2015)

\bibitem[{Locatelli(2016)}]{Locatelli2016}
Locatelli, M.: Exactness conditions for an sdp relaxation of the extended trust
  region problem. Optimization Letters \textbf{10}(6), 1141--1151 (2016)

\bibitem[{Megretski(2001)}]{10.1007/978-3-0348-8362-7_15}
Megretski, A.: Relaxations of quadratic programs in operator theory and system
  analysis. In: Borichev, A.A., Nikolski, N.K. (eds.) Systems, Approximation,
  Singular Integral Operators, and Related Topics, pp. 365--392, Birkh{\"a}user
  Basel, Basel (2001), ISBN 978-3-0348-8362-7

\bibitem[{Mixon et~al.(2016)Mixon, Villar, and Ward}]{mixon2016clustering}
Mixon, D.G., Villar, S., Ward, R.: Clustering subgaussian mixtures by
  semidefinite programming. arXiv preprint arXiv:1602.06612  (2016)

\bibitem[{Modaresi and Vielma(2017)}]{modaresi2017convex}
Modaresi, S., Vielma, J.P.: Convex hull of two quadratic or a conic quadratic
  and a quadratic inequality. Mathematical Programming \textbf{164}(1-2),
  383--409 (2017)

\bibitem[{Nesterov(1997)}]{nesterov1997quality}
Nesterov, Y.: Quality of semidefinite relaxation for nonconvex quadratic
  optimization. Tech. rep., Universit{\'e} catholique de Louvain, Center for
  Operations Research and~… (1997)

\bibitem[{Rujeerapaiboon et~al.(2019)Rujeerapaiboon, Schindler, Kuhn, and
  Wiesemann}]{rujeerapaiboon2019size}
Rujeerapaiboon, N., Schindler, K., Kuhn, D., Wiesemann, W.: Size matters:
  Cardinality-constrained clustering and outlier detection via conic
  optimization. SIAM Journal on Optimization \textbf{29}(2), 1211--1239 (2019)

\bibitem[{Santana and Dey(2018)}]{santana2018convex}
Santana, A., Dey, S.S.: The convex hull of a quadratic constraint over a
  polytope. arXiv preprint arXiv:1812.10160  (2018)

\bibitem[{Sheriff(2013)}]{sheriff2013convexity}
Sheriff, J.L.: The convexity of quadratic maps and the controllability of
  coupled systems. Ph.D. thesis (2013)

\bibitem[{Shor(1990)}]{shor1990dual}
Shor, N.Z.: Dual quadratic estimates in polynomial and boolean programming.
  Annals of Operations Research \textbf{25}(1), 163--168 (1990)

\bibitem[{Sturm and Zhang(2003)}]{SturmZhang2003}
Sturm, J.F., Zhang, S.: On cones of nonnegative quadratic functions.
  Mathematics of Operations Research \textbf{28}(2), 246--267 (2003)

\bibitem[{Tawarmalani et~al.(2002)Tawarmalani, Sahinidis, and
  Sahinidis}]{tawarmalani2002convexification}
Tawarmalani, M., Sahinidis, N.V., Sahinidis, N.: Convexification and global
  optimization in continuous and mixed-integer nonlinear programming: theory,
  algorithms, software, and applications, vol.~65. Springer Science \& Business
  Media (2002)

\bibitem[{Wang and
  K{\i}l{\i}n\c{c}-Karzan(2019{\natexlab{a}})}]{wang2019generalized}
Wang, A.L., K{\i}l{\i}n\c{c}-Karzan, F.: The generalized trust region
  subproblem: solution complexity and convex hull results. Tech. Rep.
  arXiv:1907.08843, ArXiV (2019{\natexlab{a}}),
  \urlprefix\url{https://arxiv.org/abs/1907.08843}

\bibitem[{Wang and
  K{\i}l{\i}n\c{c}-Karzan(2019{\natexlab{b}})}]{wang2019tightness}
Wang, A.L., K{\i}l{\i}n\c{c}-Karzan, F.: {On the tightness of SDP relaxations
  of QCQPs}. Tech. Rep. arXiv:1911.09195, ArXiV (2019{\natexlab{b}}),
  \urlprefix\url{https://arxiv.org/abs/1911.09195}

\bibitem[{Wolkowicz et~al.(2012)Wolkowicz, Saigal, and
  Vandenberghe}]{wolkowicz2012handbook}
Wolkowicz, H., Saigal, R., Vandenberghe, L.: Handbook of semidefinite
  programming: theory, algorithms, and applications, vol.~27. Springer Science
  \& Business Media (2012)

\bibitem[{Ye(1999)}]{ye1999approximating}
Ye, Y.: Approximating quadratic programming with bound and quadratic
  constraints. Mathematical programming \textbf{84}(2), 219--226 (1999)

\bibitem[{Ye and Zhang(2003)}]{YeZhang2003}
Ye, Y., Zhang, S.: New results on quadratic minimization. SIAM Journal on
  Optimization \textbf{14}(1), 245--267 (2003)

\bibitem[{Y{\i}ld{\i}ran(2009)}]{yildiran2009convex}
Y{\i}ld{\i}ran, U.: Convex hull of two quadratic constraints is an {LMI} set.
  IMA Journal of Mathematical Control and Information \textbf{26}(4), 417--450
  (2009)

\bibitem[{Y{\i}ld{\i}z and Cornu\'ejols(2015)}]{YC15}
Y{\i}ld{\i}z, S., Cornu\'ejols, G.: Disjunctive cuts for cross-sections of the
  second-order cone. Operations Research Letters \textbf{43}(4), 432–--437
  (2015)

\end{thebibliography}

\end{document}